\documentclass[a4paper,12pt]{article}
\usepackage[T2A]{fontenc}
\usepackage[utf8]{inputenc} 
\usepackage{amssymb,amsfonts}
\usepackage{amsfonts,amsmath,amssymb,amscd,latexsym}
\usepackage{srcltx}

\linethickness{0.4pt}

\newcommand{\bi}{\bibitem}
\newcommand{\nb}{\newblock}

\newcommand{\be}[1]{\begin{equation}\label{#1}}
\newcommand{\ee}{\end{equation}}

\newcommand{\la}{\langle\,}
\newcommand{\ra}{\,\rangle}

\newtheorem{thm}{\quad Theorem}
\newtheorem{lm}{\quad Lemma}

\title{Amenability of semigroups and the Ore condition for semigroup rings}
\author{\vspace{2ex}
V. S. Guba\thanks{This work is supported by the Russian Foundation
for Basic Research, project no. 20-01-00465.}\\
Vologda State University,\\
15 Lenin Street,\\
Vologda\\
Russia\\
160600\\
E-mail: gubavs{@}vogu35.ru}
\date{}

\begin{document}

\maketitle

\begin{abstract}

Let $M$ be a cancellative monoid. It is known~\cite{Ta54} that if $M$ is left amenable then the monoid ring $K[M]$ satisfies Ore condition, that is, there exist nontrivial common right multiples for the elements of this ring. In~\cite{Don10} Donnelly shows that a partial converse to this statement is true. Namely, if the monoid $\mathbb Z^{+}[M]$ of all elements of $\mathbb Z[M]$ with positive coefficients has nonzero common right multiples, then $M$ is left amenable. He asks whether the converse is true for this particular statement.

We show that the converse is false even for the case of groups. If $M$ is a free metabelian group, then $M$ is amenable but the Ore condition fails for $\mathbb Z^{+}[M]$. Besides, we study the case of the monoid $M$ of positive elements of R.\,Thompson's group $F$. The amenability problem for it is a famous open question. It is equivalent to left amenability of the monoid $M$. We show that for this case the monoid $\mathbb Z^{+}[M]$ does not satisfy Ore condition. That is, even if $F$ is amenable, this cannot be shown using the above sufficient condition.

\end{abstract}

\section{Preliminaries}

Let us recall one of the definitions of left amenable semigroups.

Let $S$ be a semigroup. Suppose that there exists a mapping $\mu\colon{\cal P}(S)\to[0,1]$ from the power set os $S$ into the unit interval satisfying the following conditions.
\vspace{0.5ex}

1) $\mu$ is additive, that is, $\mu(A\cup B)\mu(A)+\mu(B)$ for any disjoint subsets $A,B\subseteq S$;

2) $\mu$ is left invariant, that is, $\mu(sA)=\mu(A)$ for any $s\in S$, $A\subseteq S$;

3) $\mu$ is normalized, that is, $\mu(S)=1$.
\vspace{1ex}

Another equivalent definition can be done in terms of left invariant means on the set of bounded real-valued functions on $S$ (cf.~\cite{Don10}).
\vspace{1ex}

A convenient criterion for amenability of groups was proved by F\o{}lner in~\cite{Fol}. It can be generalized to the case of cancellative semigroups~\cite{Day68}.
\vspace{0.5ex}

Let $S$ be a cancellative semigroup. It is left amenable iff there exists $\delta > 0$ such that for any nonempty finite set $A$ in $S$, there exists a nonempty finite set $E$ such that $|aE\cap E| > \delta|E|$ for any $a\in A$. 
\vspace{1ex}

One more equivalent statement says that $S$ is left amenable whenever for any nonempty finite subset $A$ in $S$ there exists an ``almost invariant'' finite set $E$. This means that for any $\varepsilon > 0$ the set $E$ can be chosen in such a way that $|aE\Delta E| < \varepsilon|E|$ for all $a\in A$.
\vspace{1ex}

The definition of right amenable semigroups is similar. For the case of groups, both properties are equivalent.
\vspace{1ex}

A semigroup $S$ has {\em common right multiplies} whenever for any $a,b\in S$ there exist $u,v\in S$ such that $au=bv$. If $S$ has a zero element (say, $S$ is a monoid ring or so), then we say that $S$ satisfies {\em Ore condition} whenever for any $a,b\in S$ there exist $u,v\in S$ such that $au=bv$ excluding the case $u=v=0$. In this case we can also say that $S$ has nonzero (or nontrivial) common right multiplies. 
\vspace{1ex}

The following well-known result belongs to Tamari~\cite{Ta54}. 
\vspace{0.5ex}

{\sl If $M$ is left amenable cancellative monoid, then for any field $K$, the monoid ring $K[M]$ satisfies Ore condition.}
\vspace{1ex}

Let $\mathbb Z^{+}[M]$ denote the set of linear combinations of elements of $M$ with positive integer coefficients. The monoid ring $\mathbb Z[M]$ here is assumed to be cancellative. In~\cite{Don10} Donnelly proves that the following partial converse to the above statement holds. 
\vspace{1ex}

{\sl If $\mathbb Z^{+}[M]$ has nonzero common right multiples, then $M$ is left amenable.}
\vspace{1ex}

It is asked in~\cite{Don10} whether the converse to this particular statement is true. Also he askes whether the Ore condition for $K[M]$, where $K$ is a field, implies that $\mathbb Z^{+}[M]$ satisfies Ore condition. We answer both questions in negative in the next Section (even for the case of groups).
\vspace{1ex}

We are also going to discuss the above conditions for the case of R.\,Thompson's group $F$ and its positive monoid $M$. Let us recall some definitions.

Let $M$ be a monoid given by the following monoid presentation:

\be{xinf}
\la x_0,x_1,x_2,\ldots\mid x_j{x_i}=x_ix_{j+1}\ (0\le i < j)\,\ra.
\ee

It is well known that $M$ is cancellative and has common right multiplies. By Ore theorem, it is embeddable into a group given by the same presentation. This group usually denoted by $F$ was found by Richard J. Thompson in the 60s. We refer to the
survey \cite{CFP} for details. (See also \cite{BS,Bro,BG}.) This is the group of right quotients of $M$ so that $F=MM^{-1}$. We refer to $M$ as the positive monoid of $F$.

The group $F$ has presentation with 2 generators and 2 defining relations. Brin and Squier proved in~\cite{BS} that there are no free subgroups of rank $\ge2$ in $F$. The famous open problem whether $F$ is amenable, usually attributed to Geoghegan, exists for more than 40 years (see~\cite{Ger87}). It is equivalent to the left amenability of the monoid $M$ for which is known it is not right amenable~\cite{Gri90}.

It follows from the Tamari result that amenability of $F$ implies that the group ring $K[F]$ over any field satisfies Ore condition. Recently Kielak~\cite{Ba19} showed that the converse is also true. From elemenary reasons it follows that the Ore condition for the group ring $K[F]$ is equivalent to the Ore condition for monoid ring $K[M]$. So this is an open question, and we are going to clarify whether $\mathbb Z^{+}[M]$ satisfies Ore condition, where $M$ is the positive monoid of $F$. In this case amenability of $F$ would follow from Donnelly's result. However, we will show in the next Section that Ore condition does not hold for $\mathbb Z^{+}[M]$.

\section{Main Results}

\begin{lm}
\label{apm1bpm1}
Let $M$ be a monoid embeddable into a group $G$. Suppose that the monoid $\mathbb Z^{+}[M]$ has nonzero common right multiples. Then for any $a,b\in M$ there exists a relation of the form $a^{\pm1}b^{\pm1}...a^{\pm1}b^{\pm1}=1$ that holds in $G$.
\end{lm}

{\bf Proof.} Consider the equation $(1+a)u=(1+b)v$ in $\mathbb Z^{+}[M]$. It has a nontrivial solution. Therefore, $u$ and $v$ can be presented as sums of elements of $M$. That is, there exist $g_1,...,g_m\in M$ and $h_1,...,h_n\in M$ such that
$(1+a)(g_1+\cdots+g_m)=(1+b)(h_1+\cdots+h_n)$. Notice that elements in the lists may have repetitions.

Since $g_1+\cdots+g_m+ag_1+\cdots+ag_m=h_1+\cdots+h_n+bh_1+\cdots+bh_n$, the multisets $\{g_1,...,g_m,ag_1,...,ag_m\}$ and $\{h_1,...,h_n,bh_1,...,bh_n\}$ coincide. In particular, $m=n$.

Let us consider a directed graph with the vertex set $V$ indexed by the elements of the above multiset. For any $1\le i\le n$ let us add a directed edge from $ag_i$ to $g_i$ labelled by $a$. Also add a directed edge labelled by $b$ from $bh_i$ to $h_i$ for any $i$.

We get a directed labelled graph $\Gamma$. Every vertex $v$ has exactly one outcoming egde labelled by $a^{\pm1}$ and exactly one outcoming egde labelled by $b^{\pm1}$. Therefore, the edges of $\Gamma$ form a union of disjoint cycles. The label of each of these cycles form a relation that holds in the group $G$.

The proof is complete.

\begin{lm}
\label{frmet}
Let $G$ be a free metabelian group with basis $\{a,b\}$. Then $G$ has no relations between $a$ and $b$ of the form $a^{\pm1}b^{\pm1}...a^{\pm1}b^{\pm1}=1$.
\end{lm}

{\bf Proof.} Let $R$ be a normal subgroup in the free groups $\mathbb F_m$ with $m$ generators. To any word $w$ in these generators one can assign a unique path $p=p(w)$ in the Cayley graph of $\mathbb F_m/R$ starting at the identity. The word $w$ belongs to the derived subgroup $R'$ of $R$ if and only if for any edge $e$ of the Cayley graph, the number of occurrences of $e$ in $p=p(w)$ equals the number of occurrences of $e^{-1}$ in $p$. This solves the word problem in $\mathbb F_m/R'$ provided it is solvable for $\mathbb F_m/R$. The proof can be found in~\cite{DLS}; see also~\cite[Lemma 3]{Gu08}.

Let us apply this for the case $R=\mathbb F_m'$. Suppose that the word $w=a^{\pm1}b^{\pm1}...a^{\pm1}b^{\pm1}$ represents the identity in the free metabelian group, that is, in $\mathbb F_m/R'$. Then the path $p=p(w)$ in the Cayley graph of free Abelian group $\mathbb F_m/R$ satisfies the property from the previous paragraph. 

Let $e$ be an edge from this path. Without loss of generality, assume it has label $a$. The path $p=p(w)$ has an occurrence of $e^{-1}$. Up to a cyclic shift, let $p=eqe^{-1}s$ for some paths $q$, $s$. We see that $q$ is a loop, and the label of it has the form $b^{\pm1}...a^{\pm1}b^{\pm1}$. However, the length of this word is odd so it cannot represent an identity in the free Abelian group.

The proof is complete.
\vspace{1ex}

Now we can show that the converse to Donnelly's implication from~\cite{Don10} does not hold.

\begin{thm}
\label{thm1}
There exists a left amenable cancellative monoid $M$ (actually, an amenable group) such that the monoid $\mathbb Z^{+}[M]$ does not satisfy Ore condition.
\end{thm}

{\bf Proof.} Let $M$ be the free metabelian groups on 2 generators. It is well known that all soluble groups are amenable~\cite{GrL}. Therefore, $M$ is left amenable cancellative monoid. If $\mathbb Z^{+}[M]$ satisfies Ore condition then by Lemma~\ref{apm1bpm1} the group has a relation of the form $a^{\pm1}b^{\pm1}...a^{\pm1}b^{\pm1}=1$ between its free generators. This contradicts Lemma~\ref{frmet}.
\vspace{1ex}
	
Our examples answer one more question from~\cite{Don10}. If we take free metabelian group as the monoid $M$, then for any field $K$, the group ring $K[M]$ satisfies Ore condition by the result of Tamari. Howewer, this does not imply that $\mathbb Z^{+}[M]$ has nonzero right common multiples.
\vspace{1ex}

Now we want to clarify the situation for the case of the Ore condition for $\mathbb Z^{+}[M]$, where $M$ is the positive monoid of R.\,Thompson's group $F=MM^{-1}$.

\begin{lm}
\label{altrel}
Let $x_0$, $x_1$, ... , $x_m$, ... be the standard generating set for R.\,Thompson's group $F$. Then any word of the form $w=x_{i_1}^{\pm1}x_{j_1}^{\pm1}...x_{i_k}^{\pm1}x_{j_k}^{\pm1}$ $(k\ge1)$ does not represent the identity element of $F$ provided all $i_1$, ... , $i_k$ are even and all $j_1$, ... , $j_k$ are odd.
\end{lm}

{\bf Proof.} Let us call such relations {\em alternating}, where odd end even subscripts of the generators alternate. We proceed by induction on the length of the word $w$. Let $\alpha$ be the minimal subscript on letters in $w$. It is well known that $x_i\mapsto x_{i+\alpha}$ ($i\ge0$) induces a monomorphism from $F$ into itself. Therefore one can substract $\alpha$ from all indices getting a relation $w'$ in $F$. Clearly, it is also alternating. Since the algebraic sum of exponents on $x_0$ in a relation must be zero, we obtain that $w'$ has a cyclic subword of the form $...x_0^{-1}vx_0...$, where no $x_0^{\pm1}$ occur in $v$.

Applying defining relations of the form $x_i^{-1}x_jx_i=x_{j+1}$ ($i < j$) in $F$, we see that conjugation of $v=v(x_1,x_2,...)$ by $x_0$ increases all indices of letters in $v$ by $1$. The first and the last subscripts of letters in $v$ were odd; after conjugation they will be even. The subscript on a letter before $x_0^{-1}$ stays odd; the same for the letter after $x_0$. So the resulting word has the same structure where odd and even subscripts alternate. Its length decreases so the inductive assumption can be applied. This completes the proof.

\begin{thm}
\label{thm2}
Let $M$ be positive monoid of R.\,Thomspon's group $F$. Then $\mathbb Z^{+}[M]$ does not satisfy Ore condition.
\end{thm}

{\bf Proof.} Lemma~\ref{altrel} implies that $F$ has no relations of the form $a^{\pm1}b^{\pm1}...a^{\pm1}b^{\pm1}=1$, where $a=x_0$, $b=x_1$. Therefore, by Lemma~\ref{apm1bpm1}, the equation $(1+x_0)u=(1+x_1)v$ does not have nonzero solutions in $\mathbb Z^{+}[M]$.
\vspace{1ex}

Notice that there are many soltions of the equation $(1\pm x_0)u=(1\pm x_1)v$ in the monoid ring $\mathbb Z[M]$. We can give a precise description of all their solutions. Details will appear elsewhere.
\vspace{2ex}

The author thanks R.\,I.\,Grigorchuk for useful discussions.

\end{document}